\documentclass{amsart}

\usepackage{amsmath}
\usepackage{amssymb}
\usepackage{amsfonts}
\usepackage{amsbsy}
\usepackage{amscd}
\usepackage{eucal}
\usepackage{graphicx}
\usepackage{mathrsfs}

\numberwithin{equation}{section}

\newcommand{\beq}{\begin{equation}}
\newcommand{\eeq}{\end{equation}}

\newtheorem{theorem}{Theorem}[section]

\newtheorem{corollary}[theorem]{Corollary}

\newtheorem*{theoremA}{Theorem A}

\newcommand{\rmd}{\mathrm{d}}
\newcommand{\rmi}{\mathrm{i}}
\newcommand{\Real}{\mathop{\mathrm{Re}}}

\DeclareMathOperator{\sgn}{sgn}

\newcommand{\N}{\mathbb{N}}

\newcommand{\Z}{\mathbb{Z}}
\newcommand{\C}{\mathbb{C}}

\newcommand{\tfJ}{\widetilde{\mathfrak{J}}}

\newcommand{\sC}{\mathscr{C}}

\newcommand{\sds}{\strut\displaystyle}
\newcommand{\ud}{\frac{1}{2}}
\newcommand{\uds}{{\textstyle\frac{1}{2}}}

\newcommand{\skd}{\vspace*{0.2cm}}
\newcommand{\skt}{\vspace*{0.3cm}}

\begin{document}

\title[Cylinder Functions]
{On the representation of cylinder functions}

\author[Enrico De Micheli]{Enrico De Micheli}
\address{Consiglio Nazionale delle Ricerche, Via De Marini, 6 - 16149 Genova, Italy}
\email{enrico.demicheli@cnr.it}



\begin{abstract}
In this paper, we present a mixed-type integral-sum representation of the cylinder functions 
$\sC_\mu(z)$, which holds for unrestricted complex values of the order $\mu$ and 
for any complex value of the variable $z$. Particular cases of these representations and some applications,
which include the discussion of limiting forms and representations of related functions, are also
discussed.
\end{abstract}

\maketitle

\section{Introduction}
\label{se:intro}

Cylinder functions $\sC_\mu(z)$ are solutions of the Bessel differential equation
\beq
z^2\frac{\rmd^2\sC_\mu}{\rmd z^2}+z\frac{\rmd\sC_\mu}{\rmd z}+(z^2-\mu^2)\sC_\mu = 0,
\label{I.1}
\eeq
where $\mu$ is a fixed complex number. Standard cylinder functions include the Bessel function
of the first kind $J_\mu(z)$, defined by \cite[Eq. 8, p. 40]{Watson}
\beq
J_\mu(z) \doteq \left(\frac{z}{2}\right)^\mu
\sum_{k=0}^\infty \frac{(-1)^k\,(z/2)^{2k}}{k!\,\Gamma(\mu+k+1)}
\qquad (z\in\C,\mu\in\C,\mu\neq -1,-2,\ldots),
\label{I.2}
\eeq
the series on the r.h.s. of \eqref{I.2} being 
convergent absolutely and uniformly on any compact domain of $z\in\C$ and in any bounded domain 
of $\mu$. The function $J_\mu(z)$ is therefore an analytic function of $z$, except for the branch
point $z=0$ if $\mu$ is not an integer. The couple of functions $\{J_\mu(z),J_{-\mu}(z)\}$ are linearly 
independent and their linear combination gives a general solution to \eqref{I.1} if $\mu\not\in\Z$. 
When $\mu$ is an integer one is led to introduce the Bessel functions of the second kind $Y_\mu(z)$ 
(also known as Neumann's or Weber's functions), defined by \cite[Eq. 10.2.3]{DLMF}
\begin{align}
& Y_\mu(z) \doteq \frac{J_\mu(z)\cos\mu\pi-J_{-\mu}(z)}{\sin\mu\pi} & \mathrm{if}\quad\mu\not\in\Z, 
\label{I.3} \\
& Y_m(z)\doteq\lim_{\mu\to m} Y_\mu(z) & \mathrm{if}\quad m\in\Z. \label{I.4}
\end{align}
$\{J_\mu(z), Y_\mu(z)\}$ constitutes a linearly
independent pair of solutions to equation \eqref{I.1} for arbitrary $\mu\in\C$, and therefore 
the general solution to \eqref{I.1} can be written as
$\sC_\mu(z)=c_1 J_\mu(z) + c_2 Y_\mu(z)$ ($c_1, c_2$ constants). Among these combinations, 
an important role is played by the Bessel functions of the third kind (also known as Hankel's functions)
$H^{(1)}_\mu(z)$ and $H^{(2)}_\mu(z)$, in view of the asymptotic behavior for large $z$, 
which results to be very useful in applications. They are defined by \cite[Eq. 10.4.3]{DLMF}
\beq
H^{(1)}_\mu(z) \doteq J_\mu(z) + \rmi Y_\mu(z) \quad\mathrm{and}\quad 
H^{(2)}_\mu(z) \doteq J_\mu(z) - \rmi Y_\mu(z).
\label{I.5}
\eeq
In the case of purely imaginary values of $z$, ther general solution to \eqref{I.1} can be obtained 
from the functions \cite[Eqs. 10.27.6, 10.27.8]{DLMF}
\beq
I_\mu(z)\doteq e^{-\rmi\mu\pi/2}J_\mu(\rmi z) \quad\mathrm{and}\quad
K_\mu(z)\doteq
\begin{cases}
\frac{\pi}{2}\rmi^{\mu+1} H_\mu^{(1)}(\rmi z) & \mathrm{if~} -\pi<\arg z \leqslant \frac{\pi}{2}, \\[+5pt]
\frac{\pi}{2}(-\rmi)^{\mu+1} H_\mu^{(2)}(-\rmi z) & \mathrm{if~} -\frac{\pi}{2} <\arg z \leqslant \pi,
\end{cases}
\label{I.5.a}
\eeq
which are usually referred to as \emph{modified Bessel functions}. The reader is referred to the classical
monograph of G.N. Watson \cite{Watson} for more information and details on Bessel functions.

\skt

Cylinder functions enjoy several integral representations (see, e.g., \cite[Chapter 6]{Watson}),
either in terms of definite integrals and of contour integrals.
Among the representations of the Bessel function of the first kind $J_\mu(z)$ given by integrals 
on the real line (see, e.g., refs. \cite[Chapter 10.9]{DLMF}, \cite[Sect. 7.3]{Erdelyi} and 
\cite[Sect. 8.41]{Gradshteyn}), most of them hold only for restricted values of $\mu$. 
For instance, Poisson-type integral representations \cite[Eq. 7.3(3)]{Erdelyi} hold for 
$\Real\mu>-\ud$, Gubler's representations \cite[Eq. 7.3(11)]{Erdelyi} for $\Real\mu<\ud$, 
Mehler-Sonine integral formulae hold for $|\Real\mu|<\ud$ \cite[Eq. 7.12(12)]{Erdelyi} or for 
$|\Real\mu|<1$ \cite[Eq. 7.12(14)]{Erdelyi}. Only the well-known Schl\"afli's and Heine's 
representations hold for unrestricted complex values of the index $\mu$. Schl\"afli's 
representation \cite[Eq. 7.3(9)]{Erdelyi},
\beq
J_\mu(z) =\frac{1}{\pi}\int_0^\pi\cos(z\sin t-\mu t)\,\rmd t-\frac{\sin\mu\pi}{\pi}
\int_0^\infty e^{-(z\sinh t+\mu t)}\,\rmd t,
\label{Schlafli}
\eeq
holds for $\mu\in\C$, with $\Real z>0$ (it holds also in case $\Real z=0$ provided that $\Real\mu>0$),
while the two Heine's expressions \cite[Eqs. 7.3(31) \& (32)]{Erdelyi} hold for $\mu\in\C$, but separately
in the upper and lower half-plane of the $z$-plane, respectively. Both representations, Schl\"afli's 
and Heine's, are made up of the sum of two integrals.
From these representations, one can also derive corresponding integral representations of Bessel functions 
of the second and third kinds by using the relations \eqref{I.3} and \eqref{I.5}. 

In Ref. \cite{DeMicheli} we presented the following integral representation of Fourier-type for 
the Bessel functions of the first kind $J_\mu(z)$ with order $\mu\in\C$, limited to the 
half-plane $\Real\mu>-\ud$:

\begin{theoremA}[{\cite[Theorem 1]{DeMicheli}}]
\label{the:Jp}
Let $\ell\in\N_0\doteq\{0,1,2,\ldots\}$ and $\nu$ be any complex number such that $\Real\nu>-\ud$.
Then, the following integral representation for the Bessel functions of the first kind $J_{\nu+\ell}(z)$ holds:
\beq
J_{\nu+\ell}(z) = (-\rmi)^\ell \int_{-\pi}^\pi \mathfrak{J}^{(\nu)}_z(\theta) \ e^{\rmi\ell\theta}\,\rmd\theta 
\qquad (\ell\in\N_0,\Real\nu>-\uds),
\label{J.1}
\eeq
where the $2\pi$-periodic function $\mathfrak{J}^{(\nu)}_z(\theta)$ is given by
\beq
\mathfrak{J}^{(\nu)}_z(\theta) = \frac{\rmi^\nu}{2\pi} e^{\rmi\nu[\theta-\pi\sgn(\theta)]}e^{\rmi z\cos\theta} 
\ P(\nu,-\rmi z\,(1-\cos\theta)),
\label{J.2}
\eeq
and $\sgn(\cdot)$ is the sign function, $P(\nu,w)\doteq\gamma(\nu,w)/\Gamma(\nu)$ denotes 
the normalised incomplete gamma function, 
$\gamma(\nu,w)$ being the lower incomplete gamma function.
\end{theoremA}
Representation \eqref{J.1} can be put in a more appealing form as follows.
Any complex number $\mu$ can be represented uniquely as $\mu=\langle\Real\mu\rangle + \{\mu\}$, 
namely, as the sum of \emph{integral and fractional parts}, which are defined as follows: 
$\langle\Real\mu\rangle$ is the 
\emph{nearest integer} of $\Real\mu$ (sometimes also improperly referred to as the \emph{round function} 
of $\Real\mu$), with the prescription that, for $n\in\Z$: 
\beq
\langle n+\uds\rangle = n \qquad (n\in\Z).
\label{1.2.1}
\eeq
Accordingly, the (complex) \emph{fractional part} of $\mu$ is defined as: 
$\{\mu\}\doteq\mu-\langle\Real\mu\rangle$. Note that this definition implies that in general $\{\mu\}\in\C$
with $-\ud<\Real\{\mu\}\leqslant\ud$. 

Let us now consider representation \eqref{J.1} and put $\mu=\ell+\nu$ with $\ell\in\N_0$ and $\Real\nu>-\ud$. 
Without loss of generality, we can assume $\langle\Real\nu\rangle=0$. Therefore, 
$\ell=\langle\Real\mu\rangle$, with the constraint 
$\langle\Real\mu\rangle \geqslant 0$, and $\nu=\{\mu\}$ is the (complex) fractional part of $\mu$,
with $\Real\{\mu\}>-\ud$. Hence, we have $\Real\mu>-\ud$.
Thus, with simple calculations formula \eqref{J.1} can be re-written as:
\beq
J_\mu(z) = \frac{\rmi^\mu}{\pi} 
\int_0^\pi e^{-\rmi z\cos\theta} \, P(\{\mu\},-\rmi z(1+\cos\theta))\,\cos\mu\theta \,\rmd\theta
\qquad (\Real\mu>-\uds).
\label{J.3}
\eeq
Formula \eqref{J.3} generalizes to complex values of $\mu$ ($\Real\mu>-\ud$) the classical Bessel integral
\cite[Eq. 10.9.2]{DLMF}
\beq
J_n(z) = \frac{\rmi^n}{\pi}\int_0^\pi e^{-\rmi z\cos\theta}\cos n\theta\,\rmd\theta
\qquad (n\in\Z),
\label{classical_integral}
\eeq
which holds only for integer values of $n$. The ingredient which allows for this generalization is indeed the 
(normalised) incomplete gamma function (recall that $P(0,w)=1$, $w\in\C$), whose strong connection with Bessel 
functions has long been known (see, e.g., \cite{Gautschi,Tricomi1,Tricomi2}).
The constraint $\Real\mu>-\ud$, however, prevent us to extend the representation of type \eqref{J.3} to other
cylinder functions, since these latter always involve values of $J_\mu(z)$ with $\Real\mu\leqslant-\ud$ 
(see \eqref{I.3} and \eqref{I.5}).

In this paper, we aim indeed at extending the $\mu$-range of validity of formula \eqref{J.3} 
and finding a new representation of $J_\mu(z)$ which holds true for unrestricted values of $\mu\in\C$. 
This representation, which is given is Section \ref{se:J} (see formula \eqref{rep2.0}), is no longer only 
integral but it is of \emph{mixed-type}, an integral plus a finite linear combination of inverse power of $z$. 
The first part continues to be exactly the integral
given in \eqref{J.3}, whose validity is extented to any $\mu\in\C$, while the \emph{finite sum} part is actually
different from zero only for $\Real\mu\leqslant-\ud$. From this representation, particular values of $\mu$,
limiting forms and representation of related functions, i.e. error function and Dawson's integral, are analyzed.
The representation of $J_\mu(z)$ for any $\mu\in\C$ then makes it possible to obtain a similar representation
of mixed-type for the Neumann functions $Y_\mu(z)$, $\mu\in\C$. This analysis is given in Section \ref{se:Y}.
Finally, the results for the Hankel functions $H^{(1,2)}_\mu(z)$ are given in Section \ref{se:H}.

\section{Mixed-type integral-sum representation of Bessel functions
of the first kind $J_\mu(z)$}
\label{se:J}

From \eqref{J.1} we see that, for fixed $z$ and for $\ell\geqslant 0$, the function $\rmi^\ell J_{\nu+\ell}(z)$ 
coincides with the $\ell$th Fourier coefficient of the $2\pi$-periodic function 
$\mathfrak{J}^{(\nu)}_z(\theta)$. We are then brought to consider the trigonometrical series
\beq
\mathfrak{J}^{(\nu)}_z(\theta)=\frac{1}{2\pi}\sum_{\ell=-\infty}^\infty \tfJ_\ell^{(\nu)}(z) \ e^{-\rmi\ell\theta},
\label{inv.1}
\eeq
where $\tfJ_\ell^{(\nu)}(z)$ denote, for fixed $z$, 
the $\ell$th Fourier coefficients of $\mathfrak{J}^{(\nu)}_z(\theta)$:
\beq
\tfJ_\ell^{(\nu)}(z) = \int_{-\pi}^\pi \mathfrak{J}^{(\nu)}_z(\theta) \, e^{\rmi\ell\theta}\,\rmd\theta 
\qquad (\ell\in\Z).
\label{inv.2}
\eeq
In general, representation \eqref{inv.1} cannot be written explicitly since
the Fourier coefficients with $\ell<0$ are unknown.
However, this can be done in two particular cases by exploiting the parity property of 
$\mathfrak{J}^{(\nu)}_z(\theta)$.
From equation \eqref{J.2} we see that (for fixed $z$) $\mathfrak{J}^{(\nu)}_z(\theta)$ enjoys the symmetry:
$\mathfrak{J}^{(\nu)}_z(-\theta) = e^{-\rmi 2\nu[\theta-\pi\sgn(\theta)]}\,\mathfrak{J}^{(\nu)}_z(\theta)$,
which, substituted in \eqref{inv.2}, gives
\beq
\tfJ_\ell^{(\nu)}(z) = \int_{-\pi}^\pi \mathfrak{J}^{(\nu)}_z(\theta) \, e^{\rmi 2\nu\pi\sgn(\theta)}
\,e^{-\rmi(\ell+2\nu)\theta}\,\rmd\theta.
\label{inv.4}
\eeq
Formula \eqref{inv.4} induces thus a $\ell$-index symmetry on the Fourier coefficients $\tfJ_\ell^{(\nu)}(z)$ 
only if $2\nu$ is integer. Precisely, when $\nu\equiv n$ is a nonnegative integral number the symmetry formula reads:
\beq
\tfJ_{\ell}^{(n)}(z) = \tfJ_{-\ell-2n}^{(n)}(z) \qquad (\ell\in\Z; n\in\N_0),
\label{inv.5}
\eeq
while, when $\nu$ is a nonnegative half-integer: $\nu=n+\ud$, $n\in\N_0$, we have 
\beq
\tfJ_{\ell}^{(n+\ud)}(z) = -\tfJ_{-\ell-2n-1}^{(n+\ud)}(z) \qquad (\ell\in\Z; n\in\N_0).
\label{inv.6}
\eeq
For a generic complex number $\nu=n+\xi$ ($n\in\Z, \xi\in\C$) an explicit symmetry formula like 
\eqref{inv.5} and \eqref{inv.6} is no longer available. Nevertheless, the Fourier pair 
\eqref{inv.1} and \eqref{inv.2} can be exploited 
to obtain the extended representation of $J_{\mu}(z)$ we are searching for. To this end, we make use of the 
following recurrence formula for the normalised incomplete gamma function \cite[Eq. 8.8.11]{DLMF}:
\beq
P(\xi+n,w) = P(\xi,w) -w^\xi \,e^{-w} \sum_{k=0}^{n-1} \frac{w^k}{\Gamma(\xi+k+1)}
\qquad (\xi\in\C;n=0,1,2,\ldots),
\label{recurrence}
\eeq
the sum being obviously understood to be null if $n=0$. Then, the following theorem can be proved.

\begin{theorem}
\label{the:Jall}
Let $\mu$ be any complex number. Then, the following representation for the Bessel functions of the first 
kind $J_{\mu}(z)$ holds:
\beq
\begin{split}
J_{\mu}(z)
&= \frac{\rmi^{\mu}}{\pi}\int_{0}^\pi
e^{-\rmi z\cos\theta} \, P(\{\mu\},-\rmi z(1+\cos\theta))\,\cos\mu\theta \, \rmd\theta \\
&\quad +\frac{\rmi^\mu}{\sqrt{\pi}}\,e^{\rmi z}\sum_{j=1}^{-\langle\Real\mu\rangle}
\frac{\Gamma(j+\mu-\uds)}{\Gamma(j)\,\Gamma(j+2\mu)}\,(-2\rmi z)^{j+\mu-1} \qquad (\mu\in\C),
\end{split}
\label{rep2.0}
\eeq
for $z\in\C$ if $\mu\in\Z$, otherwise $z$ belongs to the slit domain $\C\setminus (-\infty,0]$,
and the sum term is understood as zero if $\langle\Real\mu\rangle\geqslant 0$, i.e., if $\Real\mu>-\ud$.
\end{theorem}

\begin{proof}
Consider the Fourier coefficient $\tfJ_\ell^{(\nu)}(z)$ (see \eqref{inv.2}) with
$\ell\in\Z$ and $\nu=n+\xi$, with $\Real\xi>-\ud$ and $n \geqslant 0$ integer.
From \eqref{inv.2} with \eqref{J.2} and using \eqref{recurrence} we have:
\beq
\begin{split}
\tfJ_\ell^{(\nu)}(z) &= \frac{\rmi^{n+\xi}}{2\pi}\int_{-\pi}^\pi e^{\rmi(n+\xi)[\theta-\pi\sgn(\theta)]}
e^{\rmi z\cos\theta} P(n+\xi,-\rmi z(1-\cos\theta))\,e^{\rmi\ell\theta}\,\rmd\theta \\
&= \frac{(-\rmi)^n\,\rmi^\xi}{2\pi}\int_{-\pi}^\pi e^{\rmi\xi[\theta-\pi\sgn(\theta)]}e^{\rmi z\cos\theta}
\left\{P(\xi,-\rmi z(1-\cos\theta)) \right. \\
&\left.\quad -[-\rmi z(1-\cos\theta)]^\xi \, e^{\rmi z(1-\cos\theta)}
\sum_{k=0}^{n-1}\frac{[-\rmi z(1-\cos\theta)]^k}{\Gamma(\xi+k+1)}\right\}\,e^{\rmi(\ell+n)\theta}\,\rmd\theta \\
&=(-\rmi)^n\frac{\rmi^\xi}{2\pi}\int_{-\pi}^\pi e^{\rmi\xi[\theta-\pi\sgn(\theta)]}e^{\rmi z\cos\theta}
P(\xi,-\rmi z(1-\cos\theta))\,e^{\rmi(\ell+n)\theta}\,\rmd\theta \\
&\quad-\frac{(-\rmi)^n}{2\pi}z^\xi e^{\rmi z} \sum_{k=0}^{n-1}\frac{(-\rmi z)^k}{\Gamma(\xi+k+1)} 
\int_{-\pi}^\pi e^{\rmi\xi[\theta-\pi\sgn(\theta)]} (1-\cos\theta)^{\xi+k}\,e^{\rmi(\ell+n)\theta}\,\rmd\theta.
\end{split}
\label{dim.1}
\eeq
Now, we assume that $n$ and $\ell$ are such that $n+\ell\geqslant 0$, 
which amounts to admitting possible negative values of $\ell$.
In view of Theorem A,
the first term on the r.h.s. of \eqref{dim.1} is proportional to $\rmi^{n+\ell}J_{\ell+n+\xi}(z)$ and
hence we can write:
\beq
\tfJ_\ell^{(\nu)}(z) = \rmi^{\ell} J_{\ell+\nu}(z) - \Xi^{(n,\xi,\ell)}(z),
\qquad (\Real\xi>-\uds, n \geqslant 0, \ell+n \geqslant 0),
\label{dim.2}
\eeq
where
\beq
\Xi^{(n,\xi,\ell)}(z) = 
\begin{cases}
\sds\frac{(-\rmi)^n}{2\pi} z^\xi e^{\rmi z} \sum_{k=0}^{n-1}\frac{(-\rmi z)^k}{\Gamma(\xi+k+1)}\, I_k^{(n,\xi,\ell)}
& \mathrm{if}\quad n \geqslant 1, \\
0 & \mathrm{if}\quad n = 0.
\end{cases}
\label{dim.3}
\eeq
with 
\beq
I_k^{(n,\xi,\ell)} \!=\! 
\int_{-\pi}^\pi \!\! e^{\rmi\xi[\theta-\pi\sgn(\theta)]} (1-\cos\theta)^{\xi+k}\,e^{\rmi(\ell+n)\theta}\,\rmd\theta,
\quad (n\geqslant 1, 0\leqslant k \leqslant n-1; \ell+n\geqslant 0).
\label{dim.4}
\eeq
The case $n=0$ is trivial since formula \eqref{dim.2} gives again the result of Theorem A. 
Then, from now on we assume $n\geqslant 1$.
With simple algebraic manipulations the integral $I_k^{(n,\xi,\ell)}$ can be written as follows with
$0\leqslant k \leqslant n-1$ and $\Real\xi>-\uds$:
\beq
\begin{split}
I_k^{(n,\xi,\ell)} &= 2 (-1)^{n+\ell}\int_{0}^\pi (1+\cos\theta)^{\xi+k}\cos(\ell+n+\xi)\theta\,\rmd\theta \\
&= \frac{2\pi (-1)^{n+\ell}}{2^{\xi+k}}\frac{\Gamma(1+2(\xi+k))}{\Gamma(1+k-\ell-n)\,\Gamma(1+k+n+\ell+2\xi)},
\end{split}
\label{dim.5}
\eeq
where we used the formula
\beq
\int_{0}^\pi (1+\cos\theta)^a \, \cos (b\,\theta) \,\rmd\theta 
= \frac{\pi}{2^a} \frac{\Gamma(1+2a)}{\Gamma(1+a-b)\,\Gamma(1+a+b)}, \qquad (\Real a >-\uds).
\label{dim.6}
\eeq
From \eqref{dim.3}, \eqref{dim.5}, using the Legendre duplication formula for the gamma function
and re-indexing the sum setting $j=k-n$, we obtain:
\beq
\Xi^{(n,\xi,\ell)}(z) = \frac{(-1)^\ell}{\sqrt{\pi}}\,(2z)^{n+\xi}\,e^{\rmi z}
\sum_{j=-n}^{-1}\frac{\Gamma(j+n+\xi+\uds)\,(-2\rmi z)^j}{\Gamma(j+1-\ell)\,\Gamma(j+1+\ell+2n+2\xi)},
\label{dim.8}
\eeq
for $n\geqslant 1$ and $\Xi^{(0,\xi,\ell)}(z)=0$.
When $n\geqslant 1$, the condition $n+\ell\geqslant 0$ means that we are considering even the case of 
Fourier coefficients with negative values of $\ell$, $-n\leqslant\ell\leqslant -1$,
which were precluded in Theorem A.
In view of the factor $\Gamma^{-1}(j+1-\ell)$, the terms of the sum in \eqref{dim.8} are different from zero only if
$j \geqslant \ell$. The latter condition, along with $\ell\geqslant-n$,  
implies that the sum in \eqref{dim.8} starts from $j=\ell$. Finally, recalling that $\nu=n+\xi$, we obtain
for values of $\nu$ such that $\langle\Real\nu\rangle\geqslant -\ell$ and $\Real\{\nu\}>-\uds$:
\beq
\Xi^{(n,\xi,\ell)}(z) = 
\begin{cases}
\sds\frac{(-1)^\ell}{\sqrt{\pi}}\,(2z)^\nu\,e^{\rmi z}
\sum_{j=\ell}^{-1}\frac{\Gamma(j+\nu+\uds)\,(-2\rmi z)^j}{\Gamma(j+1-\ell)\,\Gamma(j+1+\ell+2\nu)} &
\mathrm{if}\quad \ell\leqslant-1, \\
0 & \mathrm{if} \quad \ell \geqslant 0.
\end{cases}
\label{dim.9}
\eeq
We can now substitute \eqref{dim.9} into \eqref{dim.2} to get:
\beq
J_{\ell+\nu}(z) = (-\rmi)^\ell\,\tfJ_\ell^{(\nu)}(z)+\frac{\rmi^\ell}{\sqrt{\pi}}\,(2z)^\nu\,e^{\rmi z}
\sum_{j=\ell}^{-1}\frac{\Gamma(j+\nu+\uds)\,(-2\rmi z)^j}{\Gamma(j+1-\ell)\,\Gamma(j+1+\ell+2\nu)},
\label{dim.10}
\eeq
the second term being null for $\ell\geqslant 0$. Given $\ell \leqslant -1$, formula \eqref{dim.10} holds true 
only for $\langle\Real\nu\rangle\geqslant -\ell$ and $\Real\{\nu\}>-\ud$.
However, for any $\ell\in\Z$ the $\nu$-domain of validity of representation \eqref{dim.10} can be analytically 
extended from the half-plane $\Real\nu>-\ell-\ud$ to the half-plane $\Real\nu>-\ud$ 
if the r.h.s. of \eqref{dim.10} is proved to be an analytic function which is locally bounded in every
compact subset of $\Real\nu>-\ud$. Let see that this is the case indeed.

For what concerns the summation in \eqref{dim.10}, poles can appear from the factor $\Gamma(j+\nu+\ud)$
only if $\nu$ is a half-integer, $\nu=n+\ud$ and, consequently, when $j+n+1\leqslant 0$. 
This implies that we have simple poles when $n \leqslant -\ell-1$. However, for these values of $n$, 
poles are correspondingly present in the denominator from the term $\Gamma(j+\ell+2n+2)$, which 
has simple poles when $2\ell+2n+2 \leqslant 0$, that is, if $n\leqslant -\ell-1$, indeed. 
The ratio between these two gamma terms remains therefore finite.

Let us now consider the term $\tfJ_\ell^{(\nu)}(z)$ in the case $\ell\leqslant -1$, 
which reads (see \eqref{inv.2} and \eqref{J.2}):
\beq
\tfJ_\ell^{(\nu)}(z)=\frac{\rmi^\nu}{2\pi}\int_{-\pi}^\pi e^{\rmi\nu[\theta-\pi\sgn(\theta)]}e^{\rmi z\cos\theta} 
\ P(\nu,-\rmi z\,(1-\cos\theta)) \, e^{\rmi\ell\theta}\,\rmd\theta.
\label{dim.11}
\eeq
First, recall that $P(\nu,w)=w^\nu \gamma^*(\nu,w)$, where:
\beq
\gamma^*(\nu,w) = \frac{1}{\Gamma(\nu)}\sum_{m=0}^\infty\frac{(-w)^m}{(\nu+m)\cdot m!} \qquad
(\nu \neq 0,-1,-2,\ldots),
\label{dim.12}
\eeq
is a function entire in both $\nu$ and $w$, and $\gamma^*(-n,w)=w^n$ for $n\in\N_0$ \cite{Tricomi1}. 
Then, for $\theta\in(-\pi,\pi]$, it is easy to see these simple bounds hold true:
$\left|e^{\rmi\nu[\theta-\pi\sgn(\theta)]}\right| \leqslant e^{2\pi|\nu|}$,
$\left|e^{\rmi z\cos\theta}\right| \leqslant e^{|z|}$,
$\left|\gamma^*(\nu,-\rmi z(1-\cos\theta))\right| \leqslant \frac{e^{2|z|}}{|\Gamma(\nu+1)|}$, and
$\left|[-\rmi z(1-\cos\theta)]^\nu \right|\leqslant |z|^{|\nu|} e^{\pi|\nu|}(1-\cos\theta)^{\Real\nu}$.
Therefore, from \eqref{dim.11} it follows:
\beq
\left|\tfJ_\ell^{(\nu)}(z)\right|\leqslant
\frac{e^{3\pi|\nu|}\,e^{3|z|}\,|z|^{|\nu|}}{\pi\,|\Gamma(\nu+1)|}\,
\int_0^\pi (1-\cos\theta)^{\Real\nu}\,\rmd\theta,
\label{dim.13}
\eeq
the integral being convergent for $\Real\nu>-\ud$. Hence, formula \eqref{dim.10} defines an analytic
function of $z$ on $\C\setminus(-\infty,0]$ (on $\C$ if $\nu\in\N_0$) for $\Real\nu>-\ud$.

Finally, in \eqref{dim.10} we put $\mu=\ell+\nu$ ($\ell=\langle\Real\mu\rangle$ and $\nu=\{\mu\}$)  
and obtain the following mixed-type representation for the Bessel function of the first kind of 
any order $\mu\in\C$:
\beq
\begin{split}
J_{\mu}(z) &= \frac{\rmi^{\{\mu\}-\langle\Real\mu\rangle}}{2\pi}\int_{-\pi}^\pi
e^{\rmi\{\mu\}[\theta-\pi\sgn(\theta)]}e^{\rmi z\cos\theta}P(\{\mu\},-\rmi z(1-\cos\theta))
e^{\rmi\langle\Real\mu\rangle\theta}\,\rmd\theta \\
&\quad +\frac{\rmi^{\langle\Real\mu\rangle}}{\sqrt{\pi}}(2z)^{\{\mu\}} e^{\rmi z} \sum_{j=\langle\Real\mu\rangle}^{-1}
\frac{\Gamma(j+\{\mu\}+\uds)\,(-2\rmi z)^j}
{\Gamma(j-\langle\Real\mu\rangle+1)\,\Gamma(j+\langle\Real\mu\rangle+1+2\{\mu\})}.
\end{split}
\label{rep1}
\eeq
Rearranging the integral and re-indexing the summation, $j \to j+\langle\Real\mu\rangle-1$, 
we finally obtain formula \eqref{rep2.0}.
\end{proof}

\begin{corollary}
From \eqref{rep2.0} and the formula in \eqref{I.5.a} it is immediate to obtain the mixed-type representation
for the modified Bessel function of the first kind:
\beq
\begin{split}
I_{\mu}(z)
&= \frac{1}{\pi}\int_{0}^\pi
e^{z\cos\theta} \, P(\{\mu\},z(1+\cos\theta))\,\cos\mu\theta \, \rmd\theta \\
&\quad +\frac{1}{\sqrt{\pi}}\,e^{-z}\sum_{j=1}^{-\langle\Real\mu\rangle}
\frac{\Gamma(j+\mu-\uds)}{\Gamma(j)\,\Gamma(j+2\mu)}\,(2z)^{j+\mu-1} 
\qquad (\mu\in\C; -\pi<\arg z\leqslant\pi/2).
\end{split}
\label{modiI}
\eeq
\end{corollary}
It is worth noting that the first term in \eqref{rep2.0} coincides with the one given in \eqref{J.3}, 
which in that case was supposed to hold only for $\Real\mu>-\ud$. The second term in \eqref{rep2.0} is a finite 
linear combination of inverse powers of $z$ (in the sense that its exponent $j+\Real\mu-1<0$)
and is non-null only when $\Real\mu\leqslant-\ud$.

\skt

For later convenience, we write formula \eqref{rep2.0} as
\beq
J_\mu(z) = B_\mu(z) + \chi_\mu(z) \qquad (\mu\in\C),
\label{rep3}
\eeq
where the integral term, which is present for all values of $\mu\in\C$, is
\beq
B_\mu(z) \doteq \frac{\rmi^{\mu}}{\pi}\int_{0}^\pi
e^{-\rmi z\cos\theta} \, P(\{\mu\},-\rmi z(1+\cos\theta))\,\cos\mu\theta \, \rmd\theta \qquad (\mu\in\C),
\label{integral}
\eeq
while the \emph{corrective term}, which is non-null only when $\Real\mu\leqslant-\uds$, is given by:
\beq
\chi_\mu(z) \doteq 
\begin{cases}
\sds\frac{\rmi^\mu}{\sqrt{\pi}}\,e^{\rmi z}\!\! \sum_{j=1}^{-\langle\Real\mu\rangle}
\frac{\Gamma(j+\mu-\uds)}{\Gamma(j)\,\Gamma(j+2\mu)} \, (-2\rmi z)^{j+\mu-1}
& \mathrm{if} \quad \Real\mu\leqslant-\uds, \\
0 & \mathrm{if} \quad \Real\mu>-\uds.
\end{cases}
\label{corrective}
\eeq

As particular cases of $\chi_\mu(z)$, we note that when $\mu\equiv m\in\Z^-\doteq\{-1,-2,\ldots\}$ 
is a negative integer, we have $\chi_m(z)=0$ in view of the zero brought by the term $\Gamma^{-1}(j+2\mu)$.
Hence, $\chi_m(z)=0$ for any $m\in\Z$.
If $\mu\equiv m+\ud$ is a negative half-integer, $m\in\Z^-$, formula \eqref{corrective} 
should be understood in the limiting form in view of the ratio of the two singular terms: 
$\Gamma(j+m)/\Gamma(j+2m+1)$. Recalling that
\beq
\lim_{\stackrel{\scriptstyle\mu\to m+\ud}{m\in\Z^-}}
\frac{\Gamma(j+\mu-\uds)}{\Gamma(j+2\mu)}=2(-1)^{m+1}\frac{\Gamma(-j-2m)}{\Gamma(-j-m+1)},
\label{corrective.2}
\eeq
we have for $m\in\Z$:
\beq
\chi_{m+\ud}(z) \doteq 
\begin{cases}
\sds\frac{2}{\sqrt{\pi}}(-\rmi)^{m+\frac{3}{2}}\,e^{\rmi z}\sum_{j=1}^{-m}
\frac{\Gamma(j-m-1)}{\Gamma(j)\,\Gamma(-j-m+1)}\,(-2\rmi z)^{\ud-j} & \mathrm{if} \quad m\leqslant-1, \\
0 & \mathrm{if} \quad m\geqslant 0.
\end{cases}
\label{corrective.3}
\eeq
It is also useful for what follows to analyze the limiting forms of $\chi_\mu(z)$ for both $z\to 0,\infty$.
From \eqref{corrective} we see that near the origin $\chi_\mu(z)$ diverges as
\beq
\chi_\mu(z) \underset{z\to 0}{\sim} \frac{1}{\Gamma(1+\mu)}\left(\frac{z}{2}\right)^{\mu}
\qquad (\Real\mu\leqslant-\uds).
\label{asy_0}
\eeq
For large values of $|z|$, $\chi_\mu(z)$ is easily seen from \eqref{corrective} to vanish as:
\beq
\chi_\mu(z) \underset{z\to\infty}{\sim}\frac{\rmi^{1+\langle\Real\mu\rangle}}{\sqrt{\pi}}
\frac{\Gamma(\{\mu\}-\ud)}{\Gamma(-\langle\Real\mu\rangle)\,\Gamma(\mu+\{\mu\})}
\,(2z)^{\{\mu\}-1}\,e^{\rmi z} \qquad (\Real\mu\leqslant-\uds).
\label{asy_inf}
\eeq

\vskip 0.5cm

In order to analyze the behavior of $B_\mu(z)$ for $z\to 0$, we first
recall that the regularized incomplete gamma function can be written as: $P(\mu,w)=w^\mu\,\gamma^*(\mu,w)$
(see \eqref{dim.12}). Now, in a neighborhood of the origin $\gamma^*(\mu,w)$ is bounded by 
$\gamma^*(\mu,0)=1/\Gamma(\mu+1)$. Then, we can write:
\beq
P(\mu,w)\underset{\scriptstyle w\to 0}{\displaystyle\sim}
\begin{cases} \sds
\frac{w^\mu}{\Gamma(\mu+1)} & \mathrm{if} \quad \mu\neq -1,-2,\ldots, \\
1  & \mathrm{if} \quad \mu = -1,-2,\ldots.
\end{cases}
\label{1.8.2}
\eeq
Inserting in formula \eqref{integral} the approximation \eqref{1.8.2} along with the 
$N$th order Taylor approximation:
$e^{-\rmi z \cos\theta} \simeq \sum_{j=0}^N (-\rmi z\cos\theta)^j/j!$, we obtain:
\beq
B_\mu(z) \underset{z\to 0}{\sim} \frac{\rmi^\mu}{\pi}\frac{(-\rmi z)^{\{\mu\}}}{\Gamma(1+\{\mu\})}
\sum_{j=0}^N\frac{(-\rmi z)^j}{j!}\int_0^\pi(\cos\theta)^j\,(1+\cos\theta)^{\{\mu\}}
\cos(\langle\Real\mu\rangle+\{\mu\})\theta\,\rmd\theta,
\label{asy_B.2}
\eeq
where we have written $\mu$ in terms of its integral and fractional parts. Now, the integrals
in \eqref{asy_B.2} are null for $j<\langle\Real\mu\rangle$
(since they are proportional to $\Gamma^{-1}(j+1-\langle\Real\mu\rangle)$). Then, it must be
$N\geqslant\langle\Real\mu\rangle$ and, as first approximation (i.e., putting $N=\langle\Real\mu\rangle$) it yields:
\beq
B_\mu(z) \underset{z\to 0}{\sim} \frac{c_\mu}{\pi\Gamma(1+\langle\Real\mu\rangle)\Gamma(1+\{\mu\})}\,z^\mu,
\label{asy_B.3}
\eeq
where
\beq
c_\mu = \int_0^\pi(\cos\theta)^{\langle\Real\mu\rangle}\,(1+\cos\theta)^{\{\mu\}}
\cos\mu\theta\,\rmd\theta \qquad (\langle\Real\mu\rangle \geqslant 0).
\label{asy_B.4}
\eeq
The computation of this latter integral can be done by using formula \eqref{dim.6} and the suitable
trigonometric power formula for the term $(\cos\theta)^{\langle\Real\mu\rangle}$. 
We obtain: $c_\mu=\pi/2^\mu$.
Then, from \eqref{rep3} and using \eqref{asy_0} we have:
\beq
J_\mu(z) \underset{z\to 0}{\sim}
\begin{cases}\sds
\frac{1}{\Gamma(1+\langle\Real\mu\rangle)\Gamma(1+\{\mu\})}\left(\frac{z}{2}\right)^\mu & \mathrm{if}\quad\Real\mu > -\uds,\\[+8pt]
\sds\frac{1}{\Gamma(1+\mu)}\left(\frac{z}{2}\right)^\mu & \mathrm{if} \quad \Real\mu \leqslant -\uds
\quad (\mu\neq -1,-2,\ldots).
\end{cases}
\label{asy_B.5}
\eeq
Numerical tests show that, for $\Real\mu>-\ud$, the \emph{classical} limiting approximation of $J_\mu(z)$, as $z\to 0$,
performs better than the one given in \eqref{asy_B.5}.

\subsection{Representation of $\boldsymbol{J_{m}(z)}$ with $\boldsymbol{m\in\Z}$}
\label{subse:JInteger}

Let $\mu=m\in\Z$. Equation \eqref{rep3} gives no new information.
In fact, we have $\chi_m(z) = 0$ for all $m\in\Z$.
For what concerns $B_m(z)$ in \eqref{integral}, it yields the well-known Jacobi-Anger integral representation 
of $J_m(z)$, i.e.:
\beq
B_m(z) = \frac{\rmi^m}{\pi}\int_0^\pi e^{-\rmi z\cos\theta}\,\cos(m\theta)\,\rmd\theta \equiv J_m(z).
\label{JI}
\eeq

\vskip 0.5 cm

\subsection{Representation of the spherical Bessel function of the first kind and Dawson's integral}
\label{subse:JHelfInteger}

Let $\mu=m+\ud$, $m\in\Z$. The integral term \eqref{integral} simply reads:
\beq
B_{m+\ud}(z) = \frac{\rmi^{m+\ud}}{\pi}\int_{0}^\pi
e^{-\rmi z\cos\theta} \, P(\uds,-\rmi z(1+\cos\theta))\,\cos(m+\uds)\theta \, \rmd\theta \qquad (m\in\Z).
\label{Bhalf.1}
\eeq
This latter formula can be written differently recalling that 
$P(\uds,w^2) = \mathrm{erf}(w)$, where $\mathrm{erf}(w)$ denotes the error function
\cite[Eq. 7.2.1]{DLMF}. Then, we have:
\beq
B_{m+\ud}(z)
= \frac{\rmi^{m+\ud}}{\pi}\int_{0}^\pi
e^{-\rmi z\cos\theta} \, \mathrm{erf}(\sqrt{-2\rmi z}\cos(\theta/2)) \,\cos(m+\uds)\theta \ \rmd\theta.
\label{Bhalf.2}
\eeq
Notice that $B_{m+\ud}(z)$ is related to the $m$th Fourier coefficient $\tilde{g}_m(z)$ of the $2\pi$-periodic 
function $g_z(\theta)=\frac{1}{2\pi}\exp(-\rmi(z\cos\theta-\theta/2)\,\mathrm{erf}(\sqrt{-2\rmi z}\cos(\theta/2))$. 
Precisely, we have: $(-\rmi)^{m+\ud}B_{m+\ud}(z) = \tilde{g}_m(z)$. Therefore, we can now invert the Fourier 
representation (see \eqref{inv.1}) and, using the symmetry $\tilde{g}_m(z)=\tilde{g}_{-m-1}(z)$, 
setting $w=\sqrt{-2\rmi z}$ and recalling that $B_{m+\ud}(z)=J_{m+\ud}(z)$ for $m\geqslant 0$ (see \eqref{rep3} 
and \eqref{corrective.3}), we have thus the following representation for the error 
function\footnote{From \eqref{error.function},
representations \cite[Eq. (9.4.20) p. 57; Eq. (9.4.21) p. 58]{Luke} and \cite[Eq. 7.6.8]{DLMF} easily follow.}:
\beq
\mathrm{erf}(w\cos(\theta/2))=2\,e^{-\ud w^2\cos\theta} \ \sum_{m=0}^\infty I_{m+\ud}(w^2/2)\cos(m+\uds) \, \theta,
\label{error.function}
\eeq
where $I_\mu(z)$ denotes the modified Bessel function of the first kind.
Expression \eqref{Bhalf.2} can be written also in terms of Dawson's integral $F(z)$ \cite[Eq. 7.2.5]{DLMF}:
\beq
F(z) \doteq e^{-z^2}\int_0^z e^{t^2}\,\rmd t.
\label{dawson.1}
\eeq
We have:
\beq
B_{m+\ud}(z) = \frac{2\,\rmi^{m+3/2}\,e^{\rmi z}}{\pi\sqrt{\pi}}\int_{0}^{\pi}\!\!\!
F(-\rmi \sqrt{-2\rmi z} \cos(\theta/2))\,\cos(m+\uds)\theta \ \rmd\theta \qquad (m\in\Z).
\label{dawson.2}
\eeq
Finally, from \eqref{corrective.3} and \eqref{dawson.2}, and recalling the definition \cite[Eq. 10.47.3]{DLMF}: 
$j_m(z)=\sqrt{\ud\pi/z}\,J_{m+\ud}(z)$, we have the following representation for the spherical Bessel
function of the first kind $j_{m}(z)$:
\beq
\begin{split}
j_{m}(z) &= \sqrt{2}\rmi^{m+\frac{3}{2}}\frac{e^{\rmi z}}{\sqrt{z}}
\left[\frac{1}{\pi}\int_{0}^{\pi}F(-\rmi \sqrt{-2\rmi z} \cos(\theta/2))\,\cos(m+\uds)\theta \, \rmd\theta \right. \\
&\left. \quad + \sds\rmi\,(-1)^m \sum_{j=1}^{-m}
\frac{\Gamma(j-m-1)} {\Gamma(j)\,\Gamma(-j-m+1)}\,(-2\rmi z)^{\ud-j}\right]
\qquad (m\in\Z).
\end{split}
\label{J12.1}
\eeq
As a by-product, assuming $m\geqslant 0$ in \eqref{J12.1} we readily get the coefficient of 
the cosine transform of the function $F(w\cos(\theta/2))$:
\beq
\frac{1}{\pi}\int_0^\pi F(w\cos(\theta/2))\,\cos(m+\uds)\theta\,\rmd\theta 
= \frac{w\,e^{-w^2/2}}{2\,\rmi^m}\,j_m(-\rmi w^2/2) \qquad (m\geqslant 0),
\label{Fcosine}
\eeq
which, inverting the cosine transform, yields the following representation of Dawson's 
integral:
\beq
F(w\cos\theta) = w\,e^{-w^2/2}\sum_{m=0}^\infty(-\rmi)^m\,j_m(-\rmi w^2/2)\cos(2m+1)\theta.
\label{InvFcosine}
\eeq
A duplication formula for the Dawson integral can be obtained from \eqref{InvFcosine}. First we put
$\theta=0$ in \eqref{InvFcosine} to get
\beq
F(w) = w\,e^{-w^2/2}\sum_{m=0}^\infty(-\rmi)^m\,j_m(-\rmi w^2/2).
\label{InvFcosine.1}
\eeq
Then, expression \eqref{InvFcosine} with $\theta=\pi/3$ along with \eqref{InvFcosine.1} yields 
the duplication formula:
\beq
F(2w) = 2F(w) + 6w\,e^{-2w^2}\sum_{m=0}^\infty(-1)^m\left[j_{6m+4}(-2\rmi w^2)-\rmi\,j_{6m+1}(-2\rmi w^2)\right].
\label{InvFcosine.2}
\eeq

\subsection{Representation of the derivative $\boldsymbol{\partial J_{\mu}(z)/\partial z}$}
\label{subse:DerJ}

The representation for the derivatives of the Bessel function $J_\mu(z)$ is easily obtained from
formula \cite[Eq. 10.6.1]{DLMF}:
\beq
J_\mu^{(1)}(z)\equiv\frac{\partial J_{\mu}}{\partial z}=\ud\left[J_{\mu-1}(z)-J_{\mu+1}(z)\right].
\label{D1}
\eeq
By using \eqref{rep3} and exploiting the fact that the regularized incomplete gamma function
in \eqref{integral} depends on the fractional part of $\mu$, we can write:
\beq
\begin{split}
J_\mu^{(1)}(z) &= \frac{\rmi^{\mu-1}}{\pi}\int_0^\pi e^{-\rmi z\cos\theta} P(\{\mu\},-\rmi z(1+\cos\theta))
\cos\theta\cos\mu\theta\,\rmd\theta \\
& \quad +\ud\left[\chi_{\mu-1}(z)-\chi_{\mu+1}(z)\right].
\end{split}
\label{D2}
\eeq
The generalization to any integer order $n\geqslant 0$ is immediate by iterating formula \eqref{D1} 
and taking into account the expression of the $n$th order central differences \cite[p. 877]{Abramowitz}.
We have:
\beq
\begin{split}
J_\mu^{(n)}(z) &= \frac{\rmi^{\mu-n}}{\pi}\int_0^\pi e^{-\rmi z\cos\theta} P(\{\mu\},-\rmi z(1+\cos\theta))
(\cos\theta)^n\cos\mu\theta\,\rmd\theta \\
& \quad +\frac{1}{2^n}\sum_{j=0}^n (-1)^j\binom{n}{j}\,\chi_{\mu-n+2j}(z) \qquad (n \geqslant 0).
\end{split}
\label{D3}
\eeq
It is worth noting that, in view of definition \eqref{corrective} of $\chi_\mu(z)$,
the \emph{sum} term in \eqref{D3} is null for $\langle\Real\mu\rangle\geqslant n$.

\section{Representation of Neumann's function $Y_\mu(z)$}
\label{se:Y}

Since formula \eqref{rep2.0} holds true for any complex value of $\mu$, a similar mixed-type
representation can be written also for the Bessel function of the second type $Y_\mu(z)$ 
by means of formula \eqref{I.3}.
Therefore, from \eqref{rep3} and recalling that the term $\chi_\mu(z)$ is zero if $\Real\mu>-\ud$ (see
\eqref{corrective}), we have, for any $\mu\in\C$, the following representation of the Bessel function 
of the second kind which, even in this case, has the structure of an integral term plus a finite sum:
\beq
Y_\mu(z) = \mathcal{Y}_\mu(z) + S(\mu)\,\frac{\chi_{[-\mu\sgn(\Real\mu)]}(z)}{\sin\mu\pi}
\qquad (\mu\in\C),
\label{Y3b}
\eeq
where the integral term is given by:
\beq
\begin{split}
&\mathcal{Y}_\mu(z) \doteq \frac{B_\mu(z)\cos\mu\pi-B_{-\mu}(z)}{\sin\mu\pi}
=\frac{\rmi^\mu}{\pi\sin\mu\pi}\int_0^\pi e^{-\rmi z\cos\theta} \\
&\quad\qquad\cdot\left[\cos\mu\pi\,P(\{\mu\},-\rmi z(1+\cos\theta))-e^{-\rmi\mu\pi}\,P(-\{\mu\},-\rmi z(1+\cos\theta))\right]
\cos\mu\theta\,\rmd\theta,
\end{split}
\label{Y4}
\eeq
\beq
S(\mu) \doteq 
\begin{cases}
\cos\mu\pi  & \mathrm{if}\quad \Real\mu<0, \\
-1 & \mathrm{if}\quad \Real\mu\geqslant 0,
\end{cases}
\label{SP}
\eeq
and $\chi_{\nu}(z)$ is given in \eqref{corrective}. Note that $\chi_{[-\mu\sgn(\Real\mu)]}(z)=0$
for $|\Real\mu|<\ud$.

\subsection{Representation of $\boldsymbol{Y_m(z)}$ with $\boldsymbol{m\in\Z}$}
\label{subse:YInteger}

The Bessel function of the second kind is particularly useful when $\mu\equiv m$ is an integer since,
in this case, we have $J_{-m}(z)=(-1)^mJ_m(z)$ and $Y_m(z)$ represents the second independent 
solution to the Bessel equation \eqref{I.1}.
The mixed-type representation for $Y_m(z)$ can be obtained evaluating the limit of \eqref{Y3b} 
for $\mu\to m$. Recalling that 
$\mu = \langle\Real\mu\rangle + \{\mu\}$, this limit can be written by setting
$\langle\Real\mu\rangle = m$ and taking the limit $\{\mu\}\to 0$. Let us consider
the case $\Real\mu\geqslant 1$, the analysis when $\Real\mu\leqslant -1$ being strictly similar.
The case $\mu\to 0$ is studied separately.
We have from \eqref{Y3b}:
\beq
Y_m(z) = \lim_{\stackrel{\scriptstyle\langle\Real\mu\rangle = m\geqslant 1}{\{\mu\} \to 0}} Y_\mu(z)
=\lim_{\stackrel{\scriptstyle\langle\Real\mu\rangle = m\geqslant 1}{\{\mu\} \to 0}} \mathcal{Y}_\mu(z)
-\lim_{\stackrel{\scriptstyle\langle\Real\mu\rangle=m\geqslant 1}{\{\mu\} \to 0}}\frac{\chi_{-\mu}(z)}{\sin\mu\pi},
\label{Yint.1}
\eeq
with the latter two limits existing finite.
For what concerns the first limit on the r.h.s. of \eqref{Yint.1}, we recall that $P(0,w)=1$ and, moreover,
\beq
\lim_{a\to 0}\frac{P(a,w)-P(-a,w)}{\sin\pi a} = -\frac{2}{\pi}\,\Gamma(0,w),
\label{Yint.2}
\eeq
where $\Gamma(\nu,w)$ denotes the upper incomplete gamma function \cite[Eq. 8.2.2]{DLMF}.
Therefore, from \eqref{Y4} we obtain:
\beq
\begin{split}
& \mathcal{Y}_m(z)
=\lim_{\stackrel{\scriptstyle\langle\Real\mu\rangle = m\geqslant 1}{\{\mu\} \to 0}} \mathcal{Y}_\mu(z) \\
&= \frac{\rmi^m}{\pi}\!\!\int_0^\pi \!\! e^{-\rmi z\cos\theta}\cos m\theta \,\lim_{\{\mu\}\to 0}\!\!
\frac{P(\{\mu\},-\rmi z(1+\cos\theta))-P(-\{\mu\},-\rmi z(1+\cos\theta))}{\sin\{\mu\}\pi}\,\rmd\theta\\
&\quad+\frac{\rmi^{m+1}}{\pi}\int_0^\pi e^{-\rmi z\cos\theta}
\,\cos m\theta \lim_{\{\mu\}\to 0} P(-\{\mu\},-\rmi z(1+\cos\theta)) \, \rmd\theta \\
& =\frac{2\,\rmi^m}{\pi^2} \int_0^\pi e^{-\rmi z\cos\theta} \,
\cos m\theta \, \left[\rmi\,\frac{\pi}{2}-\Gamma(0,-\rmi z(1+\cos\theta))\right] \,\rmd\theta \qquad (m\geqslant 1).
\end{split}
\label{Yint.3}
\eeq
Recalling representation \eqref{JI} of $J_m(z)$, formula \eqref{Yint.3} can also be written as:
\beq
\mathcal{Y}_m(z)
=-\frac{2\,\rmi^m}{\pi^2} \int_0^\pi e^{-\rmi z\cos\theta} \,
\Gamma(0,-\rmi z(1+\cos\theta))\cos m\theta \,\rmd\theta + \rmi J_m(z).
\label{Yint.5}
\eeq
Concerning the second limit on the r.h.s of \eqref{Yint.1}, we have from \eqref{corrective}:
\beq
\begin{split}
& \lim_{\stackrel{\scriptstyle\langle\Real\mu\rangle = m\geqslant 1}{\{\mu\} \to 0}}
\frac{\chi_{-\mu}(z)}{\sin\mu\pi} \\
&=\frac{\rmi^{-m}}{\sqrt{\pi}}(-2\rmi z)^{-m-1}\,e^{\rmi z}
\sum_{j=1}^m (-2\rmi z)^j\,\frac{\Gamma(j-m-\ud)}{\Gamma(j)}
\lim_{\stackrel{\scriptstyle\langle\Real\mu\rangle = m\geqslant 1}{\{\mu\} \to 0}}
\frac{1}{\Gamma(j-2\mu) \, \sin\mu\pi} \\
&=\frac{2\,e^{\rmi z}}{\rmi^m\,\pi\sqrt{\pi}}\,\sum_{j=1}^m
\frac{\Gamma(\ud-j)\,\Gamma(m+j)}{\Gamma(m-j+1)}\,(2\rmi z)^{-j}
\qquad (m\geqslant 1),
\end{split}
\label{Yint.6}
\eeq
where we used
\beq
\lim_{\stackrel{\scriptstyle\langle\Real\mu\rangle = m\geqslant 1}{\{\mu\} \to 0}}
\Gamma(j-2\mu)\,\sin\mu\pi = (-1)^{m+j+1}\frac{\pi}{2\,\Gamma(2m+1-j)},
\label{Yint.7}
\eeq
and re-inxeded the sum by setting: $j\to -j+m+1$.
Finally, plugging formulae \eqref{Yint.3} and \eqref{Yint.6}
into \eqref{Yint.1} yields the following representation for 
the Bessel function of the second kind and integer order $m \geqslant 1$:
\beq
\begin{split}
Y_m(z) &= \frac{2\,\rmi^m}{\pi}\left\{\frac{1}{\pi}\int_0^\pi e^{-\rmi z\cos\theta} \,
\left[\rmi\frac{\pi}{2}-\Gamma(0,-\rmi z(1+\cos\theta))\right]\,\cos m\theta \,\rmd\theta\right.\\
&\left.\quad-\frac{(-1)^{m}}{\sqrt{\pi}}e^{\rmi z}\sum_{j=1}^m
\frac{\Gamma(\ud-j)\Gamma(j+m)}{\Gamma(m-j+1)}\,(2\rmi z)^{-j}\right\}\qquad (m\geqslant1).
\end{split}
\label{Yint.9}
\eeq
The case $m\leqslant -1$ can be treated analogously, the only difference is that in the sum term 
$m \to -m$, the integral term remaining unchanged. 
In view of \eqref{Y3b}, when $m=0$ only the integral part is non-null. Thus, we finally obtain 
the following representation for $m\in\Z$:
\beq
Y_m(z) = -\frac{2\,\rmi^m}{\pi}\left\{\frac{1}{\pi}\int_0^\pi e^{-\rmi z\cos\theta} \,
\left[\Gamma(0,-\rmi z(1+\cos\theta))-\rmi\frac{\pi}{2}\right]\,\cos m\theta \,\rmd\theta
+\sigma_m(z)\right\},
\label{Yint.10}
\eeq
where
\beq
\sigma_m(z) = \frac{(-1)^{m}}{\sqrt{\pi}}e^{\rmi z}\sum_{j=1}^{|m|}
\frac{\Gamma(\ud-j)\Gamma(j+|m|)}{\Gamma(|m|-j+1)}\,(2\rmi z)^{-j},
\label{Yint.10.1.a}
\eeq
the sum being understood to be zero when $m=0$. Explicitly, this latter case reads:
\beq
Y_0(z)=\frac{2}{\pi^2}\int_0^\pi e^{-\rmi z\cos\theta} \, 
\left[\rmi\frac{\pi}{2}-\Gamma(0,-\rmi z(1+\cos\theta))\right]\,\rmd\theta.
\label{Yint.10.1.b}
\eeq
Formula \eqref{Yint.5}, along with formulae \eqref{Yint.1} and \eqref{Yint.6}, readily yields also
a new representation for the Hankel function of the first type of integral order $H^{(1)}_m(z)$, 
$m\in\Z$ (see formula \eqref{I.5}):
\beq
\begin{split}
H_m^{(1)}(z) =
\frac{2\,\rmi^{m-1}}{\pi}\left[\frac{1}{\pi}\int_0^\pi \!\! e^{-\rmi z\cos\theta}
\,\Gamma(0,-\rmi z(1+\cos\theta))\,\cos m\theta\,\rmd\theta +\sigma_m(z) \right]\!,
\, (m\in\Z).
\end{split}
\label{Yint.12}
\eeq
The general case of Hankel functions of complex order $\mu$ will be given in 
Section \ref{se:H}.

\subsection{Representation of the spherical Bessel function of the second kind}
\label{subse:YHalfInteger}

We start from formula \eqref{Y3b} with $\mu=m+\ud$, $m\in\Z$. From \eqref{Y4} we see that
$\mathcal{Y}_{m+\ud}(z) = (-1)^{m+1}\,B_{-m-\ud}(z)$ and, using \eqref{integral}, we obtain 
(see also \eqref{dawson.2}):
\beq
\mathcal{Y}_{m+\ud}(z)
= -\frac{2\,\rmi^{m+\ud}}{\pi\sqrt{\pi}}\,e^{\rmi z}\,\int_0^{\pi} 
F(-\rmi\sqrt{-2\rmi z}\,\cos(\theta/2))\,\cos(m+\uds)\theta\,\rmd\theta.
\label{HY.0}
\eeq
For what regards the \emph{sum} term in \eqref{Y3b}, we see from \eqref{SP} that it does not contribute 
for $m\leqslant-1$ since in this case $S(m+\ud)=0$. 
Then, we limit ourselves to consider the sum term for $m\geqslant 0$. From
\eqref{corrective} and recalling that for $\mu=-m-\ud$ its integral and fractional parts
are respectively $\langle\Real\mu\rangle=-m-1$ and $\{\mu\}=\ud$, we obtain:
\beq
\begin{split}
\chi_{-m-\ud}(z) &= \frac{\rmi^{-m-\ud}}{\sqrt{\pi}}e^{\rmi z}
\sum_{j=1}^{m+1}\frac{(-2\rmi z)^{j-m-\frac{3}{2}}}{\Gamma(j)}
\lim_{\stackrel{\scriptstyle\nu\to m}{m=0,1,2,\ldots}}
\frac{\Gamma(j-\nu-1)}{\Gamma(j-2\nu-1)} \\
&=\frac{2\,\rmi^{m-\ud}}{\sqrt{\pi}} e^{\rmi z} 
\sum_{j=1}^{m+1}\frac{\Gamma(2m-j+2)}{\Gamma(j)\Gamma(m-j+2)}\,(-2\rmi z)^{j-m-\frac{3}{2}} \qquad (m\geqslant 0),
\end{split}
\label{HY.1}
\eeq
where we used the limit
\beq
\lim_{\stackrel{\scriptstyle\nu\to m}{m=0,1,2,\ldots}}
\frac{\Gamma(j-\nu-1)}{\Gamma(j-2\nu-1)} = 2\,(-1)^m\,\frac{(2m-j+1)!}{(m-j+1)!}
\qquad (1\leqslant j\leqslant m+1).
\label{HY.2}
\eeq
Now, recalling that the spherical Bessel function of the second kind is given by:
$y_m(z)=\sqrt{\pi/(2z)}Y_{m+\ud}(z)$, we can finally write:
\beq
\begin{split}
y_{m}(z) &= \sqrt{\frac{2\rmi}{z}}\,e^{\rmi z} \, \rmi^{m}
\left[-\frac{1}{\pi}\int_0^{\pi} \!\!F(-\rmi\sqrt{-2\rmi z}\,\cos(\theta/2))\,\cos(m+\uds)\theta\,\rmd\theta
\right.\\
&\left.\quad+\sds\rmi\,(-1)^m\sum_{j=0}^{m}\frac{(2m-j)!}{j!\,(m-j)!}\,(-2\rmi z)^{j-m-\frac{1}{2}}\right]
\qquad (m\in\Z).
\end{split}
\label{HY.3bis}
\eeq

\subsection{Limiting form of $\boldsymbol{Y_\mu(z)}$ near $\boldsymbol{z=0}$}
\label{subse:OriginY}

Referring to \eqref{Y3b}, let us first consider the behavior near the origin $z=0$ of 
the \emph{sum} term in formula \eqref{Y3b}. In view of \eqref{asy_0}, 
we see that this term diverges near $z=0$ as 
$O(z^{-|\Real\mu|})$. Precisely, if $\mu$ is not integer, we have:
\beq
\chi_{[-\mu\sgn(\Real\mu)]}(z)
\underset{z\to 0}{\sim}
\frac{1}{\Gamma(1-\mu\sgn(\Real\mu))}\left(\frac{2}{z}\right)^{\mu\sgn(\Real\mu)}
\qquad (\mu\in\C\setminus\Z).
\label{AS.0}
\eeq
For what concerns the integral term $\mathcal{Y}_\mu(z)$ in \eqref{Y3b}, we first recall that
the regularized incomplete gamma function can be written as: $P(\mu,w)=w^\mu\,\gamma^*(\mu,w)$
(see \eqref{dim.12}).
Since in \eqref{Y4} we have: $-\ud<\Real\{\mu\}<\ud$, for the moment we limit ourselves to consider the case 
$\mu\neq 0$. The study of the case $\mu = 0$ will be considered later. From \eqref{1.8.2} we have:
\beq
P(\{\mu\},-\rmi z(1+\cos\theta)) \underset{\scriptstyle z\to 0}{\displaystyle\sim}
\begin{cases} \sds
0 & \mathrm{if} \quad 0 < \Real\{\mu\} < \uds, \\ \sds
\frac{[-\rmi z(1+\cos\theta)]^{\{\mu\}}}{\Gamma(1+\{\mu\})} & \mathrm{if} \quad -\uds < \Real\{\mu\} < 0.
\end{cases}
\label{1.8.3}
\eeq
Analogously,
\beq
P(-\{\mu\},-\rmi z(1+\cos\theta)) \underset{\scriptstyle z\to 0}{\displaystyle\sim}
\begin{cases} \sds
\frac{[-\rmi z(1+\cos\theta)]^{-\{\mu\}}}{\Gamma(1-\{\mu\})} & \mathrm{if} \quad 0 < \Real\{\mu\} < \uds, \\
0 & \mathrm{if} \quad -\uds < \Real\{\mu\} < 0,
\end{cases}
\label{1.8.4}
\eeq
Then, we recall the following integral:
\beq
\int_0^\pi \frac{\cos\mu\theta}{(1+\cos\theta)^{\{\mu\}}}\,\rmd\theta
=\frac{2^{\{\mu\}}\,\pi\,\Gamma(1-2\{\mu\})}{\Gamma(1-\mu-\{\mu\})\,\Gamma(1+\mu-\{\mu\})},\qquad (\Real\{\mu\}<\uds).
\label{1.8.5}
\eeq
Assume for the moment $0<\Real\{\mu\}<\ud$; the case $-\ud<\Real\{\mu\}<0$ is strictly analogous.
Using \eqref{1.8.4}, \eqref{1.8.5}, and noting that $e^{-\rmi z\cos\theta}\to 1$ as 
$z\to 0$, we have from \eqref{Y4}:
\beq
\begin{split}
\mathcal{Y}_\mu(z) & \doteq \frac{\rmi^\mu}{\pi\,\sin\mu\pi}\int_0^\pi e^{-\rmi z\cos\theta}\,\cos\mu\theta \\ 
&\qquad\cdot\left[P(\{\mu\},-\rmi z(1+\cos\theta))\cos\mu\pi-e^{-\rmi\mu\pi}P(-\{\mu\},-\rmi z(1+\cos\theta))\right]\,
\rmd\theta \\
& \underset{\scriptstyle z\to 0}{\displaystyle\sim}
-\frac{\rmi^{-\mu}}{\pi\,\sin\mu\pi}\frac{(-\rmi z)^{-\{\mu\}}}{\Gamma(1-\{\mu\})}
\,\int_0^\pi \frac{\cos\mu\theta}{(1+\cos\theta)^{\{\mu\}}}\,\rmd\theta \\
& =-\frac{\rmi^{-\mu}}{\sqrt{\pi}\sin\mu\pi}\frac{\Gamma(\ud-\{\mu\})}{\Gamma(1-\mu-\{\mu\})\,\Gamma(1+\mu-\{\mu\})}
\left(\frac{\rmi}{2z}\right)^{\{\mu\}}\quad\!\!\!\!\!\! (0<\Real\{\mu\}<\uds).
\end{split}
\label{1.8.6}
\eeq
The case $-\ud<\Real\{\mu\}<0$ is treated analogously. Summarizing, 
the approximation to \eqref{Y4} in the vicinity of $z=0$ can be written for $|\Real\mu|\geqslant \ud$:
\beq
\mathcal{Y}_\mu(z) \!\underset{\scriptstyle z\to 0}{\displaystyle\sim}\!
\begin{cases}
\sds\frac{\rmi^{\mu}\cos\mu\pi}{\sqrt{\pi}\sin\mu\pi} \frac{\Gamma(\ud+\{\mu\})}{\Gamma(1-\mu+\{\mu\})\,
\Gamma(1+\mu+\{\mu\})}\left(\frac{2z}{\rmi}\right)^{\{\mu\}}\!\!, &\!\!\!\!\!\! -\uds<\Real\{\mu\}<0, \\[+15pt]
\sds\frac{-\rmi^{-\mu}}{\sqrt{\pi}\sin\mu\pi}\frac{\Gamma(\ud-\{\mu\})}{\Gamma(1-\mu-\{\mu\})\,\Gamma(1+\mu-\{\mu\})}
\left(\frac{\rmi}{2z}\right)^{\{\mu\}}, & 0<\Real\{\mu\}<\uds.
\end{cases}
\label{YY}
\eeq
Therefore, in the vicinity of $z=0$ we have: $\mathcal{Y}_\mu(z) = O(z^{-|\Real\{\mu\}|})$, which is 
negligeable with respect to the term $\Sigma_{[-\mu\sgn(\Real\mu)]}(z)$ of formula \eqref{AS.0}. 
When $|\Real\mu|<\ud$ the latter \emph{sum} term is null (see \eqref{SP}) and the
behavior near $z=0$ is governed only by the integral term \eqref{YY} with $\mu=\{\mu\}$.
Then, the approximation for $|\Real\mu|<\ud$ reads:
\beq
\mathcal{Y}_\mu(z) \underset{\scriptstyle z\to 0}{\displaystyle\sim}
\begin{cases}
\sds\frac{\cos\mu\pi}{\sin\mu\pi}\frac{1}{\Gamma(1+\mu)}\left(\frac{z}{2}\right)^\mu 
& \mathrm{if}\quad -\uds<\Real\mu<0, \\[+15pt]
\sds-\frac{1}{\sin\mu\pi}\frac{1}{\Gamma(1-\mu)} \left(\frac{z}{2}\right)^{-\mu} 
& \mathrm{if} \quad 0<\Real\mu<\uds.
\end{cases}
\label{YY2}
\eeq
Written in more compact form, we have:
\beq
Y_\mu(z) \, \underset{\scriptstyle z\to 0}{\displaystyle\sim} \,
\frac{S(\mu)}{\sin\mu\pi} \ \frac{1}{\Gamma(1-\mu\sgn(\Real\mu))}
\left(\frac{2}{z}\right)^{\mu\sgn(\Real\mu)} \quad \!\! (z,\mu\in\C;|\Real\mu|<\uds, \mu\neq 0),
\label{YY3}
\eeq
where $S(\mu)$ is given in \eqref{SP}.
Note that formula \eqref{YY3} has the same structure as formula \eqref{AS.0} which holds for any complex (not
integer) value of $\mu$. Therefore, from \eqref{Y3b}, \eqref{AS.0} and \eqref{YY3} we can extend the validity of
\eqref{YY3} to the entire complex $\mu$-plane (but integers) and write:
\beq
Y_\mu(z) \, \underset{\scriptstyle z\to 0}{\displaystyle\sim} \,
\frac{S(\mu)}{\sin\mu\pi} \ \frac{1}{\Gamma(1-\mu\sgn(\Real\mu))}\,
\left(\frac{2}{z}\right)^{\mu\sgn(\Real\mu)} \qquad (z,\mu\in\C, \mu\not\in\Z).
\label{YY3.bis}
\eeq
Finally, we use the formula $\Gamma(z)\Gamma(1-z)=\pi/\sin\pi z$ to obtain the final form
(see \cite[Eqs. 10.7.4, 10.7.5]{DLMF}:
\beq
Y_\mu(z) \, \underset{\scriptstyle z\to 0}{\displaystyle\sim} \,
\frac{\sgn(\Real\mu)\,S(\mu)}{\pi}\Gamma(\mu\sgn(\Real\mu))\left(\frac{2}{z}\right)^{\mu\sgn(\Real\mu)}
\qquad (z,\mu\in\C, \mu\not\in\Z).
\label{YY3.tris}
\eeq

\vskip 0.5 cm

The limiting behavior of $Y_\mu(z)$ when $\mu\equiv m\in\Z$ 
is easily seen by considering formula \eqref{Yint.10}.
We first observe that, for $m\neq 0$, the integral term in \eqref{Yint.10} is finite when $z\to 0$,
precisely:
\beq
\frac{2\,\rmi^m}{\pi^2}\int_0^\pi e^{-\rmi z\cos\theta} \,
\left[\rmi\frac{\pi}{2}-\Gamma(0,-\rmi z(1+\cos\theta))\right]\,\cos m\theta \,\rmd\theta
\underset{\scriptstyle z\to 0}{\longrightarrow} \,
-\frac{2\,(-\rmi)^m}{\pi \, |m|},
\label{Yint.10.1.1}
\eeq
and, therefore, the behavior for $z\sim 0$ is dominated by the \emph{sum} term. In formula \eqref{Yint.10}
we retain the term with $j=|m|$ and then, recalling that $\Gamma(\ud-|m|)\Gamma(\ud+|m|)=(-1)^m\pi$, we obtain:
\beq
Y_m(z) \, \underset{\scriptstyle z\to 0}{\displaystyle\sim} \,
-\frac{\rmi^{(m-|m|)}}{\pi}\,\Gamma(|m|)\left(\frac{2}{z}\right)^{|m|}
\qquad (z\in\C,m\in\Z,m\neq 0),
\label{YY8.3}
\eeq
which coincides indeed with formula \eqref{YY3.tris} when $\mu\to m\in\Z\setminus\{0\}$.

For the case $m=0$ we refer to representation \eqref{Yint.10.1.b}. We recall that $\Gamma(0,w)$ coincides 
with the \emph{exponential integral} $E_1(w)$ \cite[Eq. 6.2.1]{DLMF} and, moreover, it can be written as 
\cite[Eq. 6.6.2]{DLMF}:
\beq
\Gamma(0,w) = -\gamma-\ln w-\sum_{k=1}^\infty\frac{(-w)^k}{k\cdot k!},
\label{Zero.1}
\eeq
where $\gamma$ denotes the Euler-Mascheroni constant. Considering that 
$e^{-\rmi z\cos\theta}\to 1$ for $z\to 0$, we see for the first term in \eqref{Yint.10.1.b}:  
$\frac{2}{\pi^2}\int_0^\pi\rmi\frac{\pi}{2}e^{-\rmi z\cos\theta}\rmd\theta\to\rmi$ for $z\to 0$.
For what concerns the second term in \eqref{Yint.10.1.b}, we use the following definite integrals:
\beq
\int_0^\pi\ln(-\rmi z(1+\cos\theta))\,\rmd\theta = \pi\ln\left(-\rmi z/2\right),
\quad\int_0^\pi (1+\cos\theta)^k\,\rmd\theta = \frac{2^k\sqrt{\pi}\,\Gamma(k+\ud)}{\Gamma(k+1)},
\label{Zero.3}
\eeq
and formula \eqref{Zero.1} to obtain:
\beq
\begin{split}
&-\frac{2}{\pi^2}\!\int_0^\pi \!\! e^{-\rmi z\cos\theta}\Gamma(0,-\rmi z(1+\cos\theta))\,\rmd\theta
\underset{\scriptstyle z\to 0}{\displaystyle\sim}
\frac{2\gamma}{\pi}+\frac{2}{\pi}\ln\left(\frac{z}{2\rmi}\right)
+\frac{4}{\pi}\sum_{k=1}^\infty\left(\frac{\rmi z}{2}\right)^k \frac{\Gamma(2k)}{(k!)^3}.
\end{split}
\label{Zero.2}
\eeq
Finally, from \eqref{Yint.10.1.b} we obtain the following approximation near $z=0$:
\beq
Y_0(z) \, \underset{\scriptstyle z\to 0}{\displaystyle\sim} \,
\frac{2}{\pi}[\ln(z/2)+\gamma]
+\frac{4}{\pi}\sum_{k=1}^\infty\left(\frac{\rmi z}{2}\right)^{\!\!k} \frac{\Gamma(2k)}{(k!)^3} \\
\, \underset{\scriptstyle z\to 0}{\displaystyle\sim} \,
\frac{2}{\pi}\left[\ln\left(z/2\right)+\gamma\right].
\label{Zero.5}
\eeq

\section{Representation of Hankel's functions $H^{(1,2)}_\mu(z)$}
\label{se:H}

Similarly to what we have done for the Neumann function $Y_\mu(z)$, a mixed-type representation, i.e., integral plus
finite sum, can be obtained also for the Hankel functions of the first and second kind
(also known as Bessel's functions of the third kind) $H^{(1)}_\mu(z)$ and $H^{(2)}_\mu(z)$ by using 
the formulae \cite[Eqs. 10.4.7, 10.4.8]{DLMF} (see also \eqref{I.5}):
\begin{align}
H_\mu^{(1)}(z) &= \rmi \ \frac{e^{-\rmi\mu\pi}J_\mu(z)-J_{-\mu}(z)}{\sin\mu\pi}, \label{H.1} \\
H_\mu^{(2)}(z) &= -\rmi\,\frac{e^{\rmi\mu\pi}J_\mu(z)-J_{-\mu}(z)}{\sin\mu\pi}.\label{H2.1}
\end{align}
Here, we limit ourselves to present the final representations since the procedure we follow
is the same as the one presented in Section \ref{se:Y} for the Neumann functions. 
From \eqref{H.1} and \eqref{rep3} we have:
\beq
-\rmi H^{(1)}_\mu(z) = \mathcal{H}^{(1)}_\mu(z) + T_-(\mu) \,\frac{\chi_{[-\mu\sgn(\Real\mu)]}}{\sin\mu\pi}
\qquad (\mu\in\C),
\label{H3b}
\eeq
where:
\beq
T_-(\mu) \doteq 
\begin{cases}
e^{-\rmi\mu\pi}  & \mathrm{if}\quad \Real\mu<0, \\
-1 & \mathrm{if}\quad \Real\mu\geqslant 0,
\end{cases}
\label{HSP}
\eeq
the integral term being:
\beq
\begin{split}
&\mathcal{H}^{(1)}_\mu(z) \doteq\frac{e^{-\rmi\mu\pi}B_{\mu}(z)-B_{-\mu}(z)}{\sin\mu\pi}  \\
&=\!\frac{\rmi^{-\mu}}{\pi\sin\mu\pi}\!\int_0^\pi \!\! e^{-\rmi z\cos\theta}
\left[P(\{\mu\},-\rmi z(1+\cos\theta))-P(-\{\mu\},-\rmi z(1+\cos\theta))\right]
\cos\mu\theta\,\rmd\theta. 
\end{split}
\label{H4}
\eeq
When $\mu\equiv m\in\Z$, $H^{(1)}_m(z)$ can be obtained as the limit for $\mu\to m$ of formula \eqref{H4}
and keeping into account the limit \eqref{Yint.2}. We have (see also \eqref{Yint.12}):
\beq
\begin{split}
H_m^{(1)}(z) &=
\frac{2\,\rmi^{m-1}}{\pi}
\left[\frac{1}{\pi}\int_0^\pi e^{-\rmi z\cos\theta}\,\Gamma(0,-\rmi z(1+\cos\theta))\,\cos m\theta\,\rmd\theta
+\sigma_m(z)\right],
\end{split}
\label{Yint.12.bis}
\eeq
where $\sigma_m(z)$ is the \emph{finite sum} defined in \eqref{Yint.10.1.a}.

Similarly, the representation for the Hankel function $H^{(2)}_\mu(z)$ follows from
\eqref{H2.1} and \eqref{rep3}:
\beq
\rmi H^{(2)}_\mu(z) = \mathcal{H}^{(2)}_\mu(z) + T_+(\mu)\,\frac{\chi_{[-\mu\sgn(\Real\mu)]}(z)}{\sin\mu\pi}
\qquad (\mu\in\C),
\label{H2.3b}
\eeq
where:
\beq
T_+(\mu) \doteq 
\begin{cases}
e^{\rmi\mu\pi}  & \mathrm{if}\quad \Real\mu < 0, \\
-1 & \mathrm{if}\quad \Real\mu\geqslant 0,
\end{cases}
\label{HSP2}
\eeq
the integral term being:
\beq
\begin{split}
\mathcal{H}^{(2)}_\mu(z) & \doteq\frac{e^{\rmi\mu\pi}B_{\mu}(z)-B_{-\mu}(z)}{\sin\mu\pi}
=\frac{\rmi^{\mu}}{\pi\sin\mu\pi}\int_0^\pi e^{-\rmi z\cos\theta} \\
&\quad\cdot\left[e^{\rmi\mu\pi}P(\{\mu\},-\rmi z(1+\cos\theta))-e^{-\rmi\mu\pi}
P(-\{\mu\},-\rmi z(1+\cos\theta))\right]
\cos\mu\theta\,\rmd\theta. 
\end{split}
\label{H2.4}
\eeq
When $\mu\equiv m\in\Z$, we have:
\beq
\begin{split}
H_m^{(2)}(z) &=
-\frac{2\,\rmi^{m-1}}{\pi}
\left[\frac{1}{\pi}\int_0^\pi e^{-\rmi z\cos\theta}\left[\Gamma(0,-\rmi z(1+\cos\theta))-\rmi\pi\right]
\cos m\theta\,\rmd\theta+\sigma_m(z)\right].
\end{split}
\label{H2.4.bis}
\eeq

\skd

\bibliographystyle{amsplain}

\end{document}